\documentclass[11pt, reqno]{article}
\usepackage{amsfonts}
\usepackage{amsmath}
\usepackage{amsthm}
\usepackage{amssymb}
\usepackage{mathrsfs}
\numberwithin{equation}{section}
\date{}

\theoremstyle{plain}
\newtheorem{thm}{Theorem}[section]
\newtheorem{defn}[thm]{Definition}
\newtheorem{lemma}[thm]{Lemma}
\newtheorem{corollary}[thm]{Corollary}

\newtheorem{prop}[thm]{Proposition}

\newcommand{\beq}{\begin{equation}}
\newcommand{\eeq}{\end{equation}}
\newcommand{\spann}{\text{span}\,}
\def\N{\mathbb N}
\def\R{\mathbb R}
\def\Z{\mathbb Z}
\def\C{\mathbb C}
\def\F{{\cal F}}
\newcommand{\ol}{\overline}
\newcommand{\supp}{\operatorname{supp}}

\begin{document}
\title{Radial multiresolution in dimension three}
\author{Holger Rauhut and Margit R\"osler}

\maketitle

\begin{abstract} We present a construction of a wavelet-type orthonormal basis
for the space of radial $L^2$-functions in $\R^3$ via the concept of a radial 
multiresolution
analysis. The elements of the basis are obtained from a single
radial wavelet by usual dilations  and generalized translations. Hereby the generalized translation reveals the group convolution of radial functions in $\R^3$.  We provide a simple way to construct a radial scaling function and a radial wavelet
from an even classical scaling function on $\R$. Furthermore, decomposition and 
reconstruction algorithms are formulated. 
\end{abstract}


\section{Introduction} 
Standard approaches in multivariate wavelet analysis are based 
on the construction of multiresolution analyses and wavelet bases 
 from affine transformations of a finite set of basis functions, 
called multi-wavelets (see e.g. 
\cite{M}, \cite{BMM}). The 
translations are taken from a lattice subgroup $\Gamma$ of $(\R^d, +)$ and the
 dilations
are given by the integer powers of an expansive matrix which leaves $\Gamma$ 
invariant. 
The number of multi-wavelets needed to obtain a full basis of $L^2(\R^d)$ 
depends on the  determinant of the dilation matrix and is in general larger than $1$ for 
$d>1$. 
 
However, if one restricts to the analysis of radially symmetric functions
only, it suggests itself to exploit this symmetry in the construction of
corresponding wavelet transforms in order to reduce the high amount of
computational effort.  A purely radial setting would for example naturally
occur when separating variables in polar coordinates, and treating the
spherical and radial parts separately. There is a broad literature dealing
with wavelet analysis and multiresolution on spheres, 
see e.g. \cite{AV}, 
\cite{FGS} and the references therein.  In the radial case, it is not
difficult to establish a continuous wavelet analysis based on the convolution
structure of radial functions or measures instead of the usual translation in
$\R^d$. Radial convolution structures are special cases (for half-integer indices) 
of so-called Bessel-Kingman hypergroups (see \cite{BH}, \cite{Ki}),
and a continuous wavelet analysis can in fact be developed in this  general
setting,
see e.g. \cite{T}, \cite{H}, and \cite{R}. Essentially the same concept is underlying the approach
 of Epperson and Frazier 
\cite{EF}, where radial wavelet expansions in $\R^d$ are constructed which 
are based on sampling lattices with the spatial discretization 
determined by the positive zeros of related Bessel
functions of the first kind. The spatial lattice is equidistant only in the
special cases $d=1$ and $d=3$.
 This can be seen as an obstruction against a
multiscale approach to radial wavelets in arbitrary dimensions.

As to the authors knowledge, 
radial multiresolution 
analyses have in fact not been considered up to now, and 
there seems to be no general rigorous approach available for the construction of 
orthogonal radial wavelet bases in arbitrary dimension. This problem is 
 closely
related with the question to find a Poisson summation
formula compatible with the Bessel-Kingman translation, which 
still remains open. We mention that 
 the construction of \cite{EF} 
relies on the requirement that the involved radial wavelets are band-limited,
i.e. their Fourier transforms have compact support. 

In the present paper, we construct radial multiscale analyses and 
orthogonal radial wavelet bases in $\R^3$.
In dimension $3$,
the algebraic structure of the radial convolution 
allows to carry out the constructions along the same lines as in 
 the well-known Euclidean setting. Hereby, the equidistance of
the zeros of the corresponding Bessel function (which is simply a sinc-function)
is  of decisive importance. 
In order to motivate our approach and to make the intrinsic problems
towards an extension to arbitrary dimensions more visible, we start the paper
with a short account  on the continuous wavelet transform for Bessel-Kingman hypergroups, and also 
present some material about Bessel  frames. This in particular comprises
the radial case in arbitrary dimension. The continuous
 transform is, up to normalization, the same 
as in \cite{T}, while the Bessel frames are constructed in the spirit of the radial wavelet bases in \cite{EF}. In 
particular, in the spatial discretization  the zeros of associated
Bessel functions occur in a natural way, and the wavelets are band-limited.

The radial analysis
 in $\R^3$  is then also carried out in the setting
of the corresponding Bessel-Kingman hypergroup $H$ on $[0,\infty)$.
Our concept of a radial multiresolution analysis (MRA) 
in $\R^3$ is in fact that of a MRA for 
the $L^2$-space $L^2(H)$ of this hypergroup. 
The scale spaces $(V_j)_{j\in \Z}\subset L^2(H)$ are obtained by dyadic dilations from $V_0$, which in turn  is
spanned by equidistant hypergroup translates of a fixed  ``radial'' 
scaling function $\phi\in L^2(H)$. It is characterized by a two-scale relation, but in contrast to classical MRAs, $\phi$ itself is not contained  
in $V_0$, and the scale spaces are not shift-invariant (w.r.t. the hypergroup translation). Particular 
emphasis is put on the construction of orthogonal MRAs.
Here periodicity arguments similar to the classical case are needed which would not be 
available in arbitrary dimensions. From a given orthogonal MRA we then derive an 
orthogonal wavelet basis for the underlying hypergroup. 
By construction, this ``radial'' basis has  a direct interpretation as an orthogonal wavelet basis for the subspace 
of radial functions in $L^2(\R^3)$.
We also provide a concise  characterization of 
 radial scaling functions in terms of even classical 
scaling functions on $\R$. It in particular implies that in contrast to
the classical case, there do not exist any real-valued orthonormal radial
scaling functions with compact support.

The paper is organized as follows: In Section 2, we recall basic facts from
the analysis of radial functions in $\R^d$, explain the 
 corresponding radial hypergroup convolution structure,
and extend the setting to  Bessel-Kingman hypergroups of arbitrary index.
Section 3 contains a short 
account on the continuous  wavelet transform based
on the Bessel-Kingman translation, as well as the construction of Bessel frames. 
 In Section 4, radial multiresolution analyses in $\R^3$ 
and their scaling functions are introduced and discussed, while 
 Section 5 is devoted to the construction of 
orthogonal radial wavelet bases. The connection 
between radial scaling functions in $\R^3$ and classical scaling functions on $\R$ is
established in Section 6.
Finally, in Section 7 decomposition and reconstruction algorithms are discussed.

\section{Radial analysis and Bessel-Kingman hypergroups}

Suppose $F\in L^2(\R^d)$ is radial, i.e. $F(Ax)=F(x)$ a.e. for all
$A \in SO(d).$ Then there is a unique $f\in L^2(\R_+, \omega_{d/2-1})$ such that
$F(x) = f(|x|)$, where $|\,.\,|$ denotes the Euclidean norm on $\R^d$   and 
 for $\alpha \geq 0$, the  measure   $\omega_\alpha$ on  $\R_+ =[0,\infty)\,$   is
defined by
\[ d\omega_\alpha(r) = \bigl(2^\alpha\Gamma(\alpha +1)\bigr)^{-1} r^{2\alpha +1}dr.\]
Its normalization implies that  $\|F\|_2 = \|f\|_{2,\,\omega_{d/2-1}}$, where
$\|\,.\,\|_2$ is taken with
respect to the normalized Lebesgue measure $(2\pi)^{-d/2}dx$ on $\R^d$. 
On  $L^2(\R_+,\omega_\alpha)$ the Hankel transform of index $\alpha$
 is defined by
\[ \widehat f^{\,\alpha}(\lambda) = \int_0^\infty j_{\alpha}(\lambda r) f(r)
d\omega_\alpha(r) \]
with the normalized Bessel function
\[ j_\alpha(z) = \Gamma(\alpha +1) (z/2)^{-\alpha} J_\alpha(z); \quad J_\alpha(z)  =
\,
\sum_{n=0}^\infty \frac{(-1)^n
  (z/2)^{2n + \alpha}}{n!\,\Gamma(n+\alpha+1)}.\]
There is a  Plancherel Theorem for the Hankel transform, which  states  that
$f\mapsto \widehat f^\alpha$ 
establishes a self-inverse, isometric isomorphism of $L^2(\R_+,\omega_\alpha)$.
If   
 $F\in  L^2(\R^d)$ is radial with   $F(x) = f(|x|)$, then a short calculation
 shows that 
 its Plancherel transform 
\[ \F(F)(\xi) = \frac{1}{(2\pi)^{d/2}}\int_{\R^d} F(x) 
e^{-i\langle x,\xi\rangle}dx\]
is again radial with
$\,\F(F)(\xi) = \widehat f^{(d/2-1)}(|\xi|).\,$
This is due to the fact that
\begin{equation}\label{besselsphere}
j_{d/2-1}(|z|) = \int_{S^{d-1}} e^{-i\langle z,\xi\rangle}d\sigma(\xi),
\end{equation}
where $d\sigma$ denotes the spherical surface measure normalized according to 
$d\sigma(S^{d-1})=1$.
In contrast,  the usual group translates $\, x\mapsto F(x+y), \,y\in \R^d,$
will no longer be radial (apart from trivial cases). However, we observe that the
spherical means 
\[ M_rF(x):=  \int_{S^{d-1}} F(x+r\xi) d\sigma(\xi), 
\quad r\in \R_+\]
of $F$ are again radial. Moreover, $\|M_rF\|_2\leq \|F\|_2$. Thus $M_r$ induces a
norm-decreasing
linear mapping
\[T_r:\, L^2(\R_+, \omega_{d/2-1})\to  L^2(\R_+, \omega_{d/2-1}),  
\quad T_rf(|x|):= M_rF(x),\]
where $f$ and $F$ are related as above. Put $\alpha = d/2-1.$ Then 
a short calculation in polar coordinates gives
\begin{equation}\label{translation}
T_rf(s) = \,C_\alpha
\int_0^\pi f\bigl(\sqrt{r^2 + s^2 -2rs\cos \varphi}\bigr)
\sin^{2\alpha}\!\varphi\, d\varphi \quad
\,\text{ with }\,   C_\alpha =
\,\frac{\Gamma(\alpha +1)}{\Gamma\bigl(\alpha
+\frac{1}{2}\bigr)\Gamma\bigl(\frac{1}{2}\bigr)}.
\end{equation}
This defines a norm-decreasing generalized translation on $L^2(\R_+, \omega_\alpha)$
not only for $\alpha =
d/2-1$, but also for
general $\alpha \geq -1/2$. Having harmonic analysis in  mind, we are thus
lead to introduce 
a corresponding measure algebra on $\R_+$:   For $r,s \in \R_+$ we define a
probability measure 
$\delta_r*_\alpha\delta_s$ on $\R_+$ by 
\begin{equation}\label{convo}
 \delta_r*_\alpha\delta_s(f):= \,C_\alpha
\int_0^\pi f\bigl(\sqrt{r^2 + s^2 -2rs\cos \varphi}\bigr)
\sin^{2\alpha}\!\varphi\, d\varphi, \quad 
 f\in C_c(\R_+). 
\end{equation}
($C_c(\R_+)$ denotes the space of continuous, compactly supported functions on
$\R_+$).
% also \cite[11.4.]{Wat}
The convolution \eqref{convo} of point measures extends uniquely to a
bilinear, commutative, associative
and  weakly
continuous convolution on the space $M_b(\R_+)$ of regular bounded Borel
measures on $\R_+$. It is probability-preserving and  makes $M_b(\R_+)$ 
 a commutative Banach-*-algebra with respect to total variation norm, 
with neutral element $\delta_0$ and 
the mapping $\mu\mapsto \overline\mu$ as involution. 
The pair $(\R_+,*_\alpha)$ is called the Bessel-Kingman hypergroup of index
$\alpha$. We write $H_\alpha$ instead of $\R_+$ when putting emphasis on the specific
convolution structure.
Generally speaking, a hypergroup is  a locally compact 
Hausdorff space together with a weakly continuous and probability preserving 
convolution of regular bounded Borel measures generalizing the measure algebra of a
locally compact group; it
also has a unit and an involution substituting the 
group inverse. In particular, every locally compact group is also a hypergroup.
 There is a well-established harmonic analysis 
 for commutative hypergroups, which  in the  special case of $H_{d/2-1}$  reflects 
the harmonic analysis of radial functions (and measures) in $\R^d$.
We refer the reader to \cite{Jew} or \cite{BH} for a general background on
hypergroups, 
including \cite{Ki} for the Bessel-Kingman case. 
Let us mention only some aspects which
are of importance in our context: The measure $\omega_\alpha$ is a Haar measure for
$H_\alpha$, i.e. it
satisfies
\[ \int_0^\infty T_sf\, d\omega_\alpha \,=\, \int_0^\infty f d\omega_\alpha \quad 
\text{for all }\, f\in C_c(\R_+).\]
Up to a constant factor, $\omega_\alpha$ is the 
unique positive Radon measure on $\R_+$ with this property. 
Moreover,
\begin{equation}\label{adjoint}
 \int_0^\infty (T_s f) g\, d\omega_\alpha \,=\, \int_0^\infty f (T_sg)
d\omega_\alpha\quad \text{ for }\, s\in
\R_+\end{equation}
whenever both integrals exist.
The Bessel functions satisfy the product formula
 \[ \delta_r*_\alpha\delta_s(j_\alpha) \,=\, j_\alpha(r) j_\alpha(s) \quad\text{ for
all }\, r,s \in \R_+;\]
see \cite[11.4]{Wat}. For half-integers $\alpha$, this is easily deduced
 from \eqref{besselsphere}. In fact, the
 functions $r\mapsto j_\alpha(\lambda r), \, \lambda\in \R_+$ are exactly those
which are bounded and
multiplicative with respect to $*_\alpha$. 
They constitute  the  so-called dual space 
of the hypergroup $H_\alpha$.
In this way,  the Hankel transform $f\mapsto \widehat f^\alpha$ 
on $L^2(\R_+, \omega_\alpha)$ can be interpreted
as a Plancherel transform for  $H_\alpha$. For abbreviation, we  put $L^p(H_\alpha):=
L^p(\R_+, \omega_\alpha)$ and we denote by $\langle\,.\,,\,.\,\rangle $ and
$\|\,.\,\|_2$ the scalar product and
norm in $L^2(H_\alpha)$. 
It follows easily from \eqref{adjoint} that for $f\in L^2(H_\alpha)$, 
\begin{equation}\label{TFB}
\widehat{T_r f}^{\alpha}(\lambda) \,=\, j_\alpha(\lambda r) 
              \widehat{f}^\alpha(\lambda).
\end{equation} 

\section{Continuous wavelet transform  and frames for  Bessel-Kingman hypergroups}

In order to put the 
multiresolution approach in the following sections into a suitable 
framework, we continue with a  short account on 
 the continuous wavelet transform and wavelet frames for 
Bessel-Kingman hypergroups. This in particular includes a continuous radial wavelet
transform and radial wavelet
frames in arbitrary dimensions.

\subsection{The continuous wavelet transform}

The following construction is essentially the same as  in \cite{T}, only with a 
different normalization of the dilation operators
 and the resulting wavelet transform. (In contrast to \cite{T}, we 
choose dilations to be unitary, see below).
We shall therefore be brief in our presentation, and refer
the reader to \cite[Sect. 6.III]{T} for further details.

Let  $ B(L^2(H_\alpha))$
 denote the space of continuous linear operators
on $L^2(H_\alpha).$
Besides the translation operators $T_r \in B(L^2(H_\alpha))$ introduced
in Section 2,   we 
consider the dilations 
\[ D_a f(r):= \frac{1}{a^{\alpha +1}} f\bigl(\frac{r}{a}\bigr), \,\, a>0\]
which are obviously unitary  in $L^2(H_\alpha)$. Notice also that
\begin{equation}\label{Dil}
\widehat{(D_af)}^{\,\alpha}\,=\, D_{1/a}\widehat f^{\,\alpha}.
\end{equation}
We define
\[ \pi: H_\alpha\times (0,\infty) \,\longrightarrow\, B(L^2(H_\alpha)), \quad
    \pi(r,a):=  T_r D_a.\]
It is easily checked that $\pi$ is continuous with respect to 
the weak operator topology on $B(L^2(H_\alpha))$, c.f. \cite[Prop. 6.III.7]{T}.

\begin{defn}
A function $g\in L^2(H_\alpha)$ is called admissible, if
\[ C_g:= \int_0^\infty |\widehat g^{\,\alpha}(\lambda)|^2 \frac{d\lambda}{\lambda}\,
       < \, \infty.\]
\end{defn}
For abbreviation, put 
\[ d\widetilde\omega_\alpha(r,a) := \frac{1}{a^{2\alpha +3}}da \,
 d\omega_\alpha(r).\]

\noindent
The following is a reformulation of \cite[Thm 6.III.1]{T} in terms of our notation:

\begin{thm} (Plancherel Theorem) 
If $g\in L^2(H_\alpha)$ is admissible, then for all $f\in L^2(H_\alpha),$
\[\int_{H_\alpha\times (0,\infty)} |\langle f, \pi(r,a) g\rangle|^2 d\widetilde
\omega_\alpha(r,a) \,=\,
    C_g \cdot \|f\|_{2}\,.\]
\end{thm}
Polarization further implies  for admissible
 $g_1, \,g_2$ and arbitrary  $f_1,\,f_2 \in  L^2(H_\alpha)$ 
the orthogonality relation
\[\int_{H_\alpha\times (0,\infty)} \langle f_1,  \pi(r,a) g_1\rangle
   \overline{\langle f_2,  \pi(r,a) g_2\rangle}\, d\widetilde\omega_\alpha(r,a) 
      \,=\,\langle f_1,f_2\rangle\cdot 
     \int_0^\infty \overline{\widehat g_1^{\,\alpha}(\lambda)}\,
     \widehat g_2^{\,\alpha}(\lambda)\frac{d\lambda}{\lambda}\,.\]

\begin{defn}
Let $g\in L^2(H_\alpha)$ be admissible. The mapping
\[ \Phi_g: L^2(H_\alpha)\, \longrightarrow\, L^2( H_\alpha\times (0,\infty),
                \widetilde\omega_\alpha), \quad  \Psi_g f(r,a):= 
                 \langle f,\pi(r,a)g\rangle\]
is  called the wavelet transform on $H_\alpha$ 
with analyzing wavelet $g$. For $\alpha = d/2-1$, it coincides with the continuous
wavelet transform on $\R^d$
of a radial function $f$ with radial wavelet $g$, see \cite{R}.
\end{defn}

\noindent
Some types of inversion formulas for this transform 
can be found in \cite{T}. 

\subsection{Bessel frames}
Let us now turn to possible discretizations.
In order to obtain  discrete versions of the
usual wavelet transform on $\R$, it is standard to use sampling 
lattices of the type $\,\{(nba^k, a^k),\, k,n\in \Z\}\,$
 with constants $a>1, \, b>0.$
Here the discretization of the translation parameter is in accordance with
the related group 
structure of $\R$  and is therefore (in general) 
not appropriate for  radial wavelet analysis. 
 Following Epperson and Frazier \cite{EF},
we propose lattices where the discretization of the translation parameter
 involves the positive zeros   
$\,0 <  \nu_{\alpha,1}\,<\,\nu_{\alpha,2}\,<\,\ldots \)
of the Bessel function $j_\alpha$. By a result of McMahon, these are
asymptotically distributed according to
\[\nu_{\alpha,n} = \bigl(n+\frac{\alpha}{2} -
\frac{1}{4}\bigr)\,\pi\,+\,O\Bigl(\frac{1}{n}\Bigr).\]
A standard lattice in $H_\alpha\times(0,\infty)$ is given by
\[\{\,(\nu_{\alpha,n}\,ba^k\,, a^k),\> k\in \Z, \,n\in \N\,\}, \quad (a>1, \, b>0).\]
 The 
``almost orthogonal'' radial wavelet expansions of 
Epperson and Frazier \cite{EF} are based on 
this type of sampling lattice (with $a=2).$ In the following result, 
the discretization of the dilation parameter is still rather arbitrary.

\begin{thm}
Let $Q$ be a countable subset of $(0,\infty)$ and $g\in L^2(H_\alpha)$. Assume that
$\,\text{supp}\, \widehat
g^{\,\alpha}$ is 
contained in $[0,l]$ for some $l>0$, i.e. $\widehat g^{\,\alpha} = 0$ a.e. on
$(l,\infty),$ and that there exist
constants $A,B >0$ such that
\[ A \leq \,\sum_{q\in Q} |\widehat g^{\,\alpha}(q\lambda)|^2 \,\leq \,B\quad\text{
for almost all }\,
\lambda\in \R_+\,.\]
For $n\in \Z$ and $q\in Q$ define ``wavelets'' $g_{n,q}\in L^2(H_\alpha)$ by
\[ g_{n,q} := M_n^\alpha\cdot T_{r_nq} D_q(g)\,=\, \pi(r_nq,q)g,\]
where 
\[ r_n:= \frac{1}{l} \nu_{\alpha,n} \quad\text{ and }\, M_n^\alpha = \,
\frac{2^{(1-\alpha)/2}\,\nu_{\alpha,n}^{\,\alpha}}{\sqrt{\Gamma(\alpha+1)}\,
|J_{\alpha+1}(\nu_{\alpha,n})|}\,.\]
Then the set  
 $\{g_{n,q}: \,n\in \N,\, q\in Q\}\,$ is a frame for
$L^2(H_\alpha)$ with frame bounds $Al^{2\alpha +2}$ and $Bl^{2\alpha+2}.$
This means that for $f\in L^2(H_\alpha)$,
\[ Al^{2\alpha+2}\cdot\|f\|^2_{2}\,\leq\, \sum_{q\in Q}\sum_{n\in \N}
\big|\langle g_{n,q}\,,f\rangle\big|^2\,\leq\,Bl^{2\alpha
  +2}\cdot \|f\|^2_{2}\,.\]
\end{thm}

\begin{proof} The decisive point in the proof 
is the fact that the
 normalized Fourier-Bessel functions 
\[\rho_n^\alpha(\lambda):=\, M_n^\alpha j_\alpha(\nu_{\alpha,n}\lambda),\>
n\in \N\]
 form an
  orthonormal basis of the Hilbert space $X_\alpha :=
  L^2([0,1], \omega_\alpha\vert_{[0,1]})$; see e.g. Erd\'{e}lyi et al.~\cite{Erd}. 
Using \eqref{Dil} and  \eqref{TFB}, we  write
\[ \widehat g_{n,q}^{\,\alpha}(\lambda)\,=\, M_n^\alpha j_\alpha(r_n q\lambda)
D_{1/q}\,\widehat
g^{\,\alpha}(\lambda)\,=\, \rho_n^\alpha\bigl(\frac{q}{l}\lambda\bigr)
D_{1/q}\,\widehat g^{\,\alpha}(\lambda).
\]
By the Plancherel theorem for the Hankel transform,
we obtain
\begin{align}\langle g_{n,q},f\rangle\, =\,&\, \langle \widehat
g^{\,\alpha}_{n,q},\widehat f^{\,\alpha}\rangle\, = \, \int_0^\infty
\rho_n^\alpha\bigl(\frac{q}{l}\lambda\bigr)
D_{1/q}\,\widehat g^{\,\alpha}(\lambda) \,\overline{\widehat f^{\,\alpha}}\,
d\omega_\alpha\notag\\
=\, &
\int_0^\infty \rho_n^\alpha\, D_{1/l}(\widehat g^{\,\alpha})
D_{q/l}(\overline{\widehat f^{\,\alpha}})
d\omega_\alpha\,=\, \big\langle D_{1/l}(\widehat g^{\,\alpha})
 D_{q/l}\bigl(\,\overline{\widehat f^{\,\alpha}}\bigr),
\rho_n^\alpha\big\rangle_{X_\alpha}\notag
\end{align}
where for the last identity, we used  that the 
support of $D_{1/l}(\widehat g^{\,\alpha})$ is contained in $[0,1]$. Parseval's
identity for $X_\alpha$ now
yields 
\[\sum_{n\in\N} \big|\langle g_{n,q}\,,f\rangle\big|^2 = \int_0^1
\big|D_{1/l}(\widehat
g^{\,\alpha})D_{q/l}\bigl(\,\overline{\widehat
  f^{\,\alpha}}\,\bigr)\big|^2\,d\omega_\alpha \,=\, 
l^{2\alpha +2}\int_0^\infty |\widehat g^{\,\alpha}(q\lambda)|^2\,|\widehat
f^{\,\alpha}(\lambda)|^2\,d\omega_\alpha(\lambda).\]
Hence,
\[ \sum_{q\in Q}\sum_{n\in \N} \big|\langle g_{n,q}\,,f\rangle\big|^2 \,=\,
 l^{2\alpha +2}\int_0^\infty |\widehat
f^{\,\alpha}(\lambda)|^2\Bigl(\,\sum_{q\in Q}|\widehat
g^{\,\alpha}(q\lambda)|^2\Bigr)\,d\omega_\alpha(\lambda).\]
This implies the assertion.
\end{proof}

The second condition of this theorem is rather implicit. Following e.g. Bernier and
Taylor \cite{BT}, it is
possible to obtain sufficient criteria which are easier to check by introducing
the concept
of ``frame generators''. But in order to stay concise, we restrict ourselves to the
most interesting special
case of a standard lattice as defined above. Here 
$\,Q= \{ a^k,\, k\in \Z\}$ with $a>1$.

\begin{prop}
Suppose that $g\in L^2(H_\alpha)$ satisfies the following conditions:
\begin{enumerate}\itemsep=-1pt
\item[\rm{(i)}] $\widehat g^{\,\alpha}$ has compact support which is 
contained in the open interval $(0,\infty).$
\item[\rm{(ii)}]  
$\text{ess\,inf}\,\big\{|\widehat
    g^{\,\alpha}(\lambda)|: \lambda\in [a^n, a^{n+1}]\big\}\,\geq \sigma\,$ for some
$n\in \Z$ and $\sigma >0$; 
\item[\rm{(iii)}] $\,\tau:=
    \|\widehat g^{\,\alpha}\|_\infty \,<\,\infty$.
\end{enumerate}
Then there exists a constant $M>0$ such that
\begin{equation}\label{frame}
 \sigma^2 \,\leq\,\sum_{k\in \Z} |\widehat
g^{\,\alpha}(a^k\lambda)|^2\,\leq\, M\tau^2 \quad \text{for almost all
  }\,\lambda\in \R_+\,.
\end{equation}
Consequently, the set $\,\{g_{n,a^k}:\, n\in \N,\, k\in \Z\}$ is a frame for
$L^2(H_\alpha)$ with bounds
$\sigma^2l^{2\alpha+2}$ and $M\tau^2 l^{2\alpha+2}$.
\end{prop}

\begin{proof}
We use the arguments of \cite{BT}, sect.4 in a simplified form 
which is adapted to our situation.  Put $F:= [1, a)$. 
Then the intervals $a^n F= [a^n, a^{n+1}),\,n\in \Z,$ form a disjoint cover 
of $(0,\infty)$. As $T:= \text{supp}\, \widehat g^{\,\alpha}$ is compact in
$(0,\infty)$, it is covered by
finitely many of the $a^n F$.
This implies that 
\[ M:= \sup_{\lambda\in (0,\infty)} \sharp\{k\in \Z:\, \lambda \in a^{-k}T\}\]
is finite. By (iii), this gives the upper bound in \eqref{frame}.
 The lower bound follows from (ii)
together with the 
fact that the $a^n F$ cover $(0,\infty)$. 
\end{proof}

\section{Radial Multiresolution Analysis in $\R^3$}

Radial analysis in $\R^3$ corresponds to the Bessel-Kingman hypergroup 
$H_\alpha$ with $\alpha = 1/2$. 
For convenience we shall usually omit
the subscript $1/2$ and put 
\[d\omega(r) := d\omega_{1/2}(r)  = \sqrt{\frac{2}{\pi}}\, r^2 dr, \quad 
    \widehat f:= \widehat f^{\,1/2},
      \quad j(r):=\, j_{1/2}(r) \,=\, \frac{\sin r}{r}.\]
 We further write $H$ instead of $H_{1/2}$ and denote by
 $\langle\,.\,,\,.\,\rangle$ and $\|\,.\,\|_2$ 
the scalar product and norm in $L^2(H)= L^2(H,\omega)$, respectively. 
Notice that the Bessel function $j$ is even on $\R$. Hence  it is natural to assume
the Hankel transform 
\[ \widehat f(\lambda) =\int_0^\infty j(\lambda r)f(r)d\omega(r)\]
of  $f\in L^2(H)$ 
to be continued to an even function on $\R$ as well. 
We shall always do this throughout the paper.
We also  mention that by a change of variables, the generalized translation
\eqref{translation} on $H$ can be written in the simple form
\begin{equation}\label{modtrans}
 T_r f(s) \,=\, \frac{1}{2rs} \int_{|r-s|}^{r+s} f(t) t~ dt.
\end{equation}
The non-negative zeros of the  Bessel function $j$ are given by  
\[ t_k := k\pi,\, k\in \N\,,\]
and the  normalized Fourier-Bessel functions 
\begin{equation}\label{def_rhok}
 \rho_k(r):= M_k\, j(t_k r)
   \quad\text{ with }\, M_k = 2^{1/4}\pi^{5/4}k
\end{equation}
form an orthonormal basis of the Hilbert space 
    $L^2([0,1], \omega\vert_{[0,1]}).$
This is equivalent to 
the obvious fact that the functions
\[s_k(r):= \bigl(2/\pi\bigr)^{1/4}\, r\rho_k(r) \,=\, \sqrt{2}\sin(k\pi r), 
\,\, k\in \N\]
are an orthonormal basis for  $L^2[0,1]:= L^2([0,1], dr)$. 
It will be of importance in the following that the $s_k$ are $2$-periodic.

Let us come to the definition of a radial multiresolution analysis (MRA) 
for $\R^3$, i.e. for the Bessel-Kingman hypergroup $H$. It is close to the
well-known definition of Mallat \cite{Mal} for $\R$. 
For convenience,  we introduce the notation 
\[T^{(k)}:= T_{t_k}\,=\, T_{k\pi}, \quad (k\in \N).\]
If $f\in L^2(H)$ then according to
  (\ref{TFB}),
\begin{equation}\label{trans2}
(M_kT^{(k)}f)^\wedge(\lambda) \,=\, \rho_k(\lambda)\widehat{f}(\lambda)\,=\,
\bigl(\frac{\pi}{2}\bigr)^{1/4}
\,\frac{s_k(\lambda)}{\lambda}\, \widehat f(\lambda).
\end{equation}

\begin{defn}{\bf (Radial Multiresolution Analysis)} \label{MRA}
A radial MRA  for $\R^3$  is a sequence $\{V_j\}_{j\in \Z}$ of 
closed linear subspaces of $L^2(H)$ such that
\begin{itemize}\itemsep=-1pt
\item[\rm{(1)}] $V_j \subseteq V_{j+1}$ for all $j \in \Z$;
\item[\rm{(2)}] $\bigcap_{j=-\infty}^\infty V_j ~=~ \{0\}$;
\item[\rm{(3)}] $ \bigcup_{j=-\infty}^\infty V_j$ is dense in $L^2(H)$;
\item[\rm{(4)}] $f\in V_j$ if and only if $f(2\,\cdot) \in V_{j+1}$;
\item[\rm{(5)}] There exists a function $\phi \in L^2(H)$ such that 
\[B_\phi \,:=\, \{M_k T^{(k)}\phi:\, k \in \N \}\]
 is a Riesz basis of $V_0$, i.e.
$\,\spann  B_\phi$ is dense in $V_0$ and there exist constants 
$A,B>0$ such that
\[ A \|\alpha\|_2^2 ~\leq~ 
   \|\sum_{k=1}^\infty \alpha_k M_k T^{(k)}\phi\,\|_{2}^2 ~\leq~ B
      \|\alpha\|_2^2 \]
for all $\alpha = (\alpha_k)_{k\in\N} \in l^2(\N)$; here 
$\,\|\alpha\|_2 = \bigl(\sum_{k=1}^\infty |\alpha_k|^2\bigr)^{1/2}$.
\end{itemize}
\end{defn}
 The function $\phi$ in (5) is called a scaling function for the MRA $\{V_j\}$. 
We remark explicitly that in contrast to the classical case,
$\phi$ itself is not contained in $V_0$ and  $V_0$ is not shift invariant; in fact, 
 if $f\in V_0$ then 
$T^{(k)}f \notin V_0$ for all $k$. This will be shown  in Corollary \ref{shift}
 below.

Our first aim is to  determine an orthonormal basis for $V_0$ from its 
Riesz basis, i.e. a function $\phi^*\in L^2(H)$ such that 
 $B_{\phi^*}$ consititutes
 an orthonormal basis for $V_0$. For $\phi\in L^2(H)$  we define
\[P_\phi(\lambda)~:=~ \sum_{n=-\infty}^\infty |\widehat{\phi}(\lambda + 2n)|^2,\]
which is even   and $2$-periodic on $\R$.

\begin{prop}\label{Prop_Riesz} Let $\phi\in L^2(H)$ and $A,B>0$. 
 Then 
\begin{equation}\label{Riesz_prop}
A \|\alpha\|_2^2 \leq \|\sum_{k=1}^\infty 
\alpha_k M_k T^{(k)}\phi\,\|_{2}^2
\leq \,B \|\alpha\|_2^2 \quad\text{ for all }\,\alpha \in l^2(\N)
\end{equation}
if and only if 
\begin{equation}\label{Riesz_prop2}
A \leq P_\phi(\lambda) \leq \,B\qquad \mbox{for almost all } \lambda \in \R.
\end{equation}
\end{prop}

\begin{proof}
Let  $\alpha\in l^2(\N)$ be an arbitrary finite sequence. Define
\begin{equation}\label{Four_sine}
\widetilde{\alpha} := \sum_{k=1}^\infty \alpha_k s_k \, \in L^2[0,1]. 
\end{equation}
We may regard $\widetilde\alpha$ as an odd, $2$-periodic function on $\R$. 
By the Plancherel theorem for the Hankel transform and \eqref{trans2}, 
\begin{align}
  \|\sum_{k=1}^\infty \alpha_k M_k &T^{(k)}\phi\,\|_{2}^2\,=\,
 \|\sum_{k=1}^\infty \alpha_k \rho_k\, \widehat \phi\,\|_{2}^2\,=\, 
  \int_0^\infty |\sum_{k=1}^\infty \alpha_k s_k(\lambda)|^2 
  |\widehat\phi(\lambda)|^2 d\lambda \notag\\
      =&\,  \frac{1}{2}\int_{-\infty}^\infty |\widetilde\alpha(\lambda)|^2 
       |\widehat\phi(\lambda)|^2d\lambda \,=\, 
     \frac{1}{2} \int_{-1}^1 |\widetilde\alpha(\lambda)|^2 
    P_\phi(\lambda)d\lambda\,=\,
  \int_0^1  |\widetilde\alpha(\lambda)|^2 P_\phi(\lambda)d\lambda. \notag
\end{align}
The $s_k$ forming an orthonormal basis of $L^2[0,1]$, we have
$\,\|\alpha\|_2 = \|\widetilde\alpha\|_{L^2[0,1]}$. As the finite sequences
form a dense subspace of $l^2(\N),$  this implies the assertion.
\end{proof}

\noindent
With $A=B =1$ we immediately obtain 

\begin{corollary}\label{lemOrthnorm} For $\phi\in L^2(H)$
 the following 
statements are equivalent:\parskip=-1pt
\begin{itemize} \itemsep=-1pt
\item[\rm{(1)}] The set $\,B_\phi = \{M_k T^{(k)}\phi : k \in \N\}$ is orthonormal
     in $L^2(H)$;
\item[\rm{(2)}] $P_\phi = 1$ a.e.
\end{itemize}
\end{corollary}

\noindent
For $\phi\in L^2(H),$ put
\[ V_\phi := \,\overline{\spann B_\phi}\,, \]
the closure being taken in $L^2(H).$ The set $B_\phi$ is a Riesz basis of $V_\phi$ if
and only if there exist
constants $A,B >0$ such that the equivalent conditions of 
Proposition \ref{Prop_Riesz} are satisfied. This will be a standard requirement in
the sequel, and we therefore
introduce a separate notation:

\begin{defn}  A function $\phi\in L^2(H)$  satisfies condition (RB)  if 
$B_\phi$ is a Riesz basis of $V_\phi\,$. 
\end{defn}

As before, we shall often consider functions from $L^2[0,1]$ as odd,
$2$-periodic functions on $\R$. We therefore define
\[ S\,:= \,\{\alpha:\,\R \to \C\,|\,\, \alpha\vert_{[0,1]} \in L^2[0,1],~ 
\alpha(-x) = -\alpha(x),~ \alpha(x+2) = \alpha(x) \,\,\text{for almost all }\,x\}.\] 
$S$ is a Hilbert space with norm $\|\,.\,\|_{L^2[0,1]}$.

\begin{lemma}\label{Lem_Contain} Let $\phi\in L^2(H)$ satisfy (RB).
 Then for $f\in L^2(H)$ 
we have the equivalence
\[f\in V_\phi\,\,\Longleftrightarrow\, \,\widehat f(\lambda) = 
   \frac{\beta(\lambda)}{\lambda} \widehat\phi(\lambda) \,\text{ with }
    \beta\in S.\]
The function $f\in V_\phi$ corresponding to  
$\beta = \sum_{k=1}^\infty {\alpha_ks_k}\in S$ 
with $(\alpha_k)_{k\in\N}\in l^2(\N)$ is given by
$f = \bigl(\frac{2}{\pi}\bigr)^{1/4}\sum_{k=1}^\infty \alpha_k M_kT^{(k)}\phi$. 
\end{lemma}

\begin{proof} By \eqref{trans2}      
we have
\[\frac{s_k(\lambda)}{\lambda}\widehat\phi(\lambda)\,=\,
   \bigl(\frac{2}{\pi}\bigr)^{1/4}
  (M_kT^{(k)}\phi)^\wedge(\lambda).\]  
The translates $M_kT^{(k)}\phi$ form a Riesz basis of $V_\phi$, and hence  the
$(M_kT^{(k)}\phi)^\wedge\,$
are a Riesz basis of  $\widehat{V_\phi}$. As $(s_k)_{k\in \N}$ is an orthonormal
basis of $S$, this implies the
assertion.
\end{proof}

\noindent
Lemma \ref{Lem_Contain} is  of particular interest when $\phi$ is the scaling
function of a MRA $\{V_j\}$. Then
$V_0 = V_\phi$, and we easily deduce
the previously mentioned lack of shift-invariance: 

\begin{corollary}\label{shift}
Let $\{V_j\}_{j\in\Z}$ be a radial MRA. Then $f\in V_0$ implies that $\,T^{(k)}\!f= 
T_{t_k}f\notin V_0$ for all $k\in \N$. Similarly, $f\in V_j$ implies that 
$ T_{2^{-j}t_k}f\notin V_j$ for all $k\in \N$.
\end{corollary}

\begin{proof}
After rescaling it is enough to consider $V_0$. Recall that 
$\,  M_k(T^{(k)}\!f)^\wedge \,=\, \rho_k \widehat f.$ But if $\beta\in S$, then 
$\rho_k\beta\notin S$ for all $k$, because periodicity is lost. The characterization
of $V_\phi = V_0$ according
to the previous Lemma thus shows that for $f\in V_0$,  $T^{(k)}\!f \notin V_0$.
\end{proof}

\noindent
In the situation of the Lemma, we can easily determine  an orthonormal basis
of $V_\phi$ by renormalization:

\begin{thm}\label{cor_orth} (Orthogonalization) 
Suppose $\phi\in L^2(H)$ satisfies condition (RB).   
Define 
$\phi^*\in L^2(H)$ by its
Hankel transform
\begin{equation}\label{orth_scale}
\widehat{\phi^*}\,:=\, \frac{\widehat{\phi}}
{\sqrt{P_\phi}}\,.
\end{equation}
Then $\,B_{\phi^*} = \{M_k T^{(k)} \phi^*\,:\,k\in\N\}$ forms an orthonormal basis of
$V_\phi = V_{\phi^*}\,$.
\end{thm}

\noindent
If $\phi$ is a scaling function of a MRA  $\{V_j\}$, then 
$\,V_\phi = V_{\phi^*} = V_0$, and we call 
$\phi^*$   an orthonormal scaling 
function for  $\{V_j\}$. 
 
\begin{proof} By definition of $\phi^*$ we have $P_{\phi^*} = 1$ a.e. and hence 
$B_{\phi^*}$ is orthonormal according to Corollary \ref{lemOrthnorm}.
It remains to prove that $\,V_{\phi^*} = V_\phi$.
For this, we have to verify that $\,M_kT^{(k)}\phi^*\in V_\phi$ and
 $\,M_kT^{(k)}\phi\in V_{\phi^*}$ for all $k\in\N$. 
Employing Lemma \ref{Lem_Contain}, relation \eqref{trans2} and finally the relation 
$\,s_k(r) = (\pi/2)^{1/4}\,r\rho_k(r)$, 
one obtains
that the above conditions are equivalent to 
\[ \frac{s_k}{\sqrt{P_\phi}}\in S, \quad s_k\sqrt{P_\phi}\in S
    \quad\text{for all }k\in\N.\]
But these conditions 
are obviously satisfied by our assumption on $P_\phi\,.$ 
\end{proof}

Let us return  to our definition of a radial MRA for $\R^3$.
Suppose we  start with a function $\phi\in L^2(H)$ satisfying condition (RB) with
Riesz constants $A,B >0.$ 
Define corresponding scale spaces $\{V_j\}_{j\in \Z}$ by 
\[
  V_0:= V_\phi\,, \quad V_j:= D_{2^{-j}} V_0\]
where the dilation operator $D_a \in B(L^2(H))$  is defined as in Section 3.
Then in particular, the $V_j$ satisfy axiom (4) of Definition \ref{MRA}.
Put further
\begin{equation}\label{bases}
\phi_{j,k}(r) := \, D_{2^{-j}}\bigl(M_kT^{(k)}\phi\bigr)(r)\,=\,8^{j/2} M_k
(T^{(k)}\phi)(2^j r)
\quad j\in \Z,\, k\in \N.
\end{equation}
Then  $\langle \phi_{j,k}, \phi_{j,l}\rangle\,=\, \langle \phi_{0,k},
\phi_{0,l}\rangle\,$ 
for all $j,k,l$. Thus the $\{\phi_{j,k},\, k\in \N\}$ 
form a Riesz basis of $V_j$, with the same Riesz constants $A,B$ as for $j=0$.
In particular,
\[ V_j \,=\,\overline{\spann\{\phi_{j,k}, k\in\N \}}.\]  
Moreover,  if  $B_\phi = \{\phi_{0,k}:\,k\in\N\}$ is an orthonormal basis for $V_0$,
then 
$\,\{\,\phi_{j,k}, \, k\in \N\}\,$
is an orthonormal basis of $V_j$. 

Recall now axiom (1)  of Def. \ref{MRA}, which 
requires that the  $V_j$ are nested. Similar as in the classical case,
this condition can be reformulated
 in terms of a two-scale relation for $\phi$:

\begin{prop}\label{Prop_Contain} 
For $\phi$ with (RB) and $\{V_j\}_{j\in \Z}$ as above, the following statements
are equivalent: \parskip=-2pt
\begin{enumerate} \itemsep=-2pt
\item[\rm{(1)}] $V_j\subseteq V_{j+1}$ for all $j\in \Z$
\item[\rm{(2)}] $V_{-1}\subseteq V_0$
\item[\rm{(3)}] There exists a function $\gamma\in S$ such that 
\begin{equation}\label{two_scale}
   \sin(2\pi\lambda) \widehat\phi(2\lambda) \,=\, 
       \gamma(\lambda)\widehat\phi(\lambda).\end{equation}
\end{enumerate}
In this case, the coefficients $(h_k)_{k\in\N}\in l^2(\N)$ in the two-scale relation
\begin{equation}\label{def_hk}
\phi_{-1,1} \,=\, \sum_{k=1}^\infty h_k\,\phi_{0,k}
\end{equation}
are the coefficients in the Fourier sine series of $\gamma \in S$:
\[\gamma \,=\,\frac{1}{2}\sum_{k=1}^\infty h_k s_k.\]
\end{prop}

\begin{proof} Rescaling by the factor $2^j$ shows that (1) and (2) are equivalent. 
For (2),  we need at least 
$\phi_{-1,1} \in V_0$. According
to Lemma \ref{Lem_Contain} this is equivalent to the existence of a 
function $\beta \in S$ such that for almost all $\lambda$, 
\begin{equation}\label{rel1}
 \sqrt{8} M_1\, j(2 \pi \lambda) 
\widehat\phi(2\lambda) \,=\, \widehat{\phi_{-1,1}}(\lambda) \,=\, 
\frac{\beta(\lambda)}{\lambda} \widehat\phi(\lambda).
\end{equation}
Moreover, if $\phi_{-1,1}$ has the expansion \eqref{def_hk}, then 
$\,\beta = (\pi/2)^{1/4}\sum_{k=1}^\infty h_ks_k$. 
In turn, equation \eqref{rel1}  is equivalent to relation \eqref{two_scale}
with  $\gamma = (8\pi)^{-1/4}\beta$. 
This gives the stated connection between $\gamma$ and $\phi_{-1,1}$. 
It remains to show that $\phi_{-1,1}\in V_0$ 
(or equivalently, relation \eqref{two_scale}) 
already implies that $\phi_{-1,k}\in V_0$ for all $k\in \N$. 
Similar as above, the latter is equivalent to 
\begin{equation}\label{rel2}
 \sin(2k\pi\lambda) \widehat\phi(2\lambda)\,=\, 
k\gamma_k(\lambda)\widehat\phi(\lambda)
\end{equation}
with $\gamma_k\in S$. The relation between $\gamma_k$ and $\phi_{-1,k}$ 
is now given by
\[\gamma_k \,=\, \frac{1}{2k}\sum_{l=1}^\infty h_l^{(k)} s_k, \qquad 
 \phi_{-1,k} \,=\, \sum_{l=1}^\infty h_l^{(k)} \phi_{0,l}.\]
Comparison of  (\ref{two_scale}) with (\ref{rel2}) yields
\begin{equation}\label{Tk_contain}
\gamma_k(\lambda)~=~ \gamma(\lambda) 
\frac{\sin(2k\pi \lambda)}{k\sin(2\pi \lambda)} ~=~
\gamma(\lambda) U_{k-1}(\cos 2 \pi \lambda) 
\end{equation}
where 
\[U_k(x) = \frac{\sin(k+1)t}{(k+1)\sin t}\,,\,\,\, x=\cos t\]
 denotes 
the $k$-th Chebychev polynomial of the second kind, normalized such that
$U_k(1)=1$. Thus 
given $\gamma\in S$, we define 
\[ \gamma_k(\lambda):= \gamma(\lambda)U_{k-1}(\cos 2\pi\lambda).\]
As $U_{k-1}$ is bounded on $[-1,1]$, $\gamma_k$ is contained in $S$ as well,
 and hence
$\phi_{-1,k}\in V_0.$ 
\end{proof}

Let us now consider the remaining axioms (2) and (3) of a radial multiresolution 
analysis.

\begin{thm}\label{thmMRA} 
Let $\phi \in L^2(H)$ satisfy  condition (RB) and assume 
 that the scale spaces
\[V_j\, = \,\overline{\spann\{\phi_{j,k}: k\in\N \}},\quad j\in \Z\]  
satisfy $\, V_{-1}\subseteq V_0$.
Suppose further that 
$\big|\widehat\phi\,\big|$ is continuous in $0$. 
Then  $\{V_j\}_{j\in \Z}$ is a radial MRA  
if and only if $\,\widehat\phi(0) \neq 0$. Moreover,
$\phi$ is an orthonormal scaling function  if and only if
$\,\big|\widehat\phi(0)\big|=1.$
\end{thm}

We remark that
continuity of $\widehat\phi$ in $0$ (even on  $\R_+$) is 
for example guaranteed if $\phi\in L^2(H)\cap L^1(H).$

\begin{proof} We have to check axioms (2) and (3). This may be done by slight 
modifications
of standard arguments in the affine case. The condition on $\widehat\phi\,$ in $0$ 
will be needed only for (3).
We define an 
orthonormal scaling function $\phi^*$ according to Theorem \ref{cor_orth}.
The orthogonal 
projection $P_j$ of $L^2(H)$ onto $V_j$ is then given by
\[ 
P_j f\,=\, \sum_{k=1}^\infty \langle f, \phi^*_{j,k} \rangle \phi^*_{j,k}\]
 where 
the $\phi^*_{j,k}$ are defined as in \eqref{bases}. 
For $(2)$, we need to show that 
$\, \lim_{j\to -\infty} \|P_j f\|_{2} = 0\,$ for all $ f\in L^2(H).\,$
 Since functions with
compact support are dense in $L^2(H),$ we may assume that 
$\text{supp} f$ is  contained in a compact interval $[0,R]$.
Parseval's equation then implies
\begin{align}
\|P_j f\|^2_{2} \,=\,& \sum_{k=1}^\infty |\langle f, \phi^*_{j,k}\rangle|^2\, \leq 
\, \sum_{k=1}^\infty \Big|\int_0^R f(r)\overline{\phi^*_{j,k}(r)}
d\omega(r)\Big|^2\notag\\
\leq\,& \|f\|_{2}^2\cdot \int_0^R \sum_{k=1}^\infty|\phi^*_{j,k}(r)|^2 d\omega(r)\,
=\,\|f\|_{2}^2\cdot\int_0^{2^j\!R} \sum_{k=1}^\infty|\phi^*_{0,k}(r)|^2
d\omega(r).\notag
\end{align}
Using the explicit formula \eqref{modtrans} for the hypergroup translation in $H$ 
we further deduce
\[ \sum_{k=1}^\infty  |\phi^*_{0,k}(r)|^2 \,=\,
\sum_{k=1}^\infty  |M_k T^{(k)}\!\phi^*(r)|^2
\,=\, \sum_{k=1}^\infty \Big|\frac{(2\pi)^{1/4}}{2r}
\int_{|k\pi -r|}^{k\pi+r} \phi^*(t)t\, dt \Big|^2 \,.\]
Now assume that $j$ is sufficiently small so that $\,2^j\!R < \pi/2$. Then for $r\in
[0,2^j\!R\,],$ the
integration domains $\,[k\pi-r,k\pi+r]\,$ do not
overlap and we obtain 
\[\sum_{k=1}^\infty  |\phi^*_{0,k}(r)|^2 \,\leq \,
\frac{C}{r}\sum_{k=1}^\infty \int_{k\pi-r}^{k\pi +r} |\phi^*(t)|^2 t^2dt\,\leq 
\frac{C^\prime}{r}\|\phi^*\|_{2}^2\]
with suitable constants $C, C^\prime >0$ independent of $j$.
Hence for $j$ sufficiently small,
\[\|P_j f\|_{2}^2 ~\leq~ C^{\prime\prime}\int_0^{2^j\!R}\frac{1}{r}
d\omega(r),\]
which tends to $0$ as $\, j\to -\infty.$ This proves (2). \hfill\break
As to  (3), suppose first that $\,\widehat\phi(0) \neq 0$ and
let $\,h \in \big(\cup_{j=-\infty}^\infty V_j\big)^\bot$, i.e.  
$P_j h = 0$ for all $j \in \Z$.
 We claim that $h=0$. Indeed, for  $\epsilon > 0$ there exists a 
function $f \in L^2(H)$ such that the support of its Hankel transform 
$\widehat f$ is compact and
$\|f - h\|_{2} \,\leq \epsilon.\,$ This implies
\[ \|P_j f\|_{2}\,=\,\|P_j(f-h)\|_{2} \,\leq \epsilon \quad\text{for all }\,j\in
\Z.\]
By the Riesz basis assumption on $\phi$, we further have
\begin{equation}\label{Riesz2}
 A \sum_{k=1}^\infty 
|\langle f, \phi_{j,k}\rangle|^2 \,\leq\, \|P_j f\|^2_{2} \,\leq\,  B
\sum_{k=1}^\infty 
|\langle f, \phi_{j,k}\rangle|^2,
\end{equation}
c.f.  Lemma 2.7 in \cite{Wojt}. Further, if $\,\text{supp}\,\widehat f\,\in [0,R]$ 
then
\[ \langle f, \phi_{j,k}\rangle\,=\, \langle\, \widehat f,\widehat \phi_{j,k} \rangle
\,=\,
\int_0^R \widehat f(\lambda) \rho_k^{(j)} (\lambda)\,
\overline{\widehat\phi(2^{-j}\lambda)}\, d\omega(\lambda)\]
where
\[\rho_k^{(j)}:=\, D_{2^j}\rho_k\,.\]
Note that the functions $\{\rho_k^{(j)}, k\in \N\}$ 
form an orthonormal basis of $L^2([0,2^j],\omega\vert_{[0,2^j]})=:X_j$. 
Suppose now
 that $j$ is sufficiently large, i.e. $2^j \geq R$. Then
\[\langle f, \phi_{j,k}\rangle\,=\, \langle\,\widehat f\,\,\overline{\widehat
\phi(2^{-j}\cdot)}\,, \rho_k^{(j)}
\rangle_{X_j}.
\]
Thus by Parseval's equation for $X_j$,
\[\sum_{k=1}^\infty |\langle f, \phi_{j,k}\rangle|^2 \,=\, 
\big\|\widehat f\,\,\overline{\widehat \phi(2^{-j}\cdot)}\big\|_{X_j}^2\,=\,
\int_0^R \big|\widehat f(\lambda)\big|^2\,
\big|{\widehat \phi(2^{-j}\lambda)}\big|^2 d\omega(\lambda). 
\]
As $\big|\widehat{\phi}\big|$ is assumed to be continuous in $0$,  the functions 
$\lambda\to\,\big|\widehat{\phi}(2^{-j} \lambda)\big|$ converge to the constant 
$|\widehat{\phi}(0)|>0\,$ uniformly on $[0,R]$ as $j \to \infty$. 
Hence 
\[\epsilon \,\geq\,  \lim_{j\to\infty} \|P_j f\|_{2}\,\geq \, \sqrt{A}\,
\big|\widehat\phi(0)\big|
\|\,\widehat f\,\|_{2} \,\geq  \, 
\sqrt{A}\,\big|\widehat{\phi}(0)\big| (\|h\|_2 - \epsilon).  \]
As $\epsilon$ is arbitrarily small, this shows that $h=0$ and hence axiom (3) 
is satisfied. Vice versa, axiom (3) implies that
\[ \lim_{j\to\infty} P_jf \,=\, f\quad\text{for all } f\in L^2(H).\]
If $\widehat f\,$ is compactly supported, then the same calculation as above shows
that 
\[ \lim_{j\to\infty} \|P_jf\|_{2} \,\leq \,\sqrt{B}\,\big|\widehat\phi(0)\big|\,
\|\widehat f\,\|_{2}\]
which enforces  $\, \widehat\phi(0) \not=0.$ Finally, the $\phi_{j,k}\,,\, k\in\N$
are
  orthonormal  if and only if  it is possible to choose $A=B=1$ in \eqref{Riesz2}. In
this case 
we obtain (for $f$ as just before)  
\[\lim_{j\to\infty} \|P_jf\|_{2} \,= \,\big|\widehat\phi(0)\big|\,
\|\widehat f\,\|_{2}\,=\, \big|\widehat\phi(0)\big|
\|f\|_{2}.\]
Thus (3) is satisfied exactly if $\big|\widehat\phi(0)\big| =1.$
\end{proof}

Let us now write the two scale relation (\ref{two_scale}) in a slightly 
different form, namely 
\begin{equation}\label{relG}
\widehat\phi(2 \lambda) \,=\, G(\lambda)\, \widehat\phi(\lambda)
\end{equation}
with 
\[G(\lambda):=\, \frac{\gamma(\lambda)}{\sin(2\pi \lambda)}.\]
The filter function $G$
is obviously $2$-periodic and even. Whenever its restriction
to $[0,1]$ is contained in $L^2[0,1]$  it can be expanded as
a cosine series,
\[G(\lambda) \,=\, \sqrt{2} \sum_{n=0}^\infty g_n \cos(n\pi \lambda).\]
As for a classical multiresolution analysis one proves the 
following.

\begin{lemma}\label{lemGprop} Suppose that $\phi \in L^2(H)$ 
is an orthonormal
scaling function of a radial MRA. 
Then  the associated filter function $G$ satisfies
\begin{equation}\label{Gsquare}
|G(\lambda)|^2 + |G(\lambda+1)|^2 \,=\, 1\quad \text{a.e.}. \end{equation}
If in addition  $\phi \in L^1(H)$ then (\ref{Gsquare}) holds
pointwise and
\[G(0)=1, \,\, G(1)= 0\]
which implies
\[ \sqrt{2}\sum_{n=0}^\infty g_n \,=\, 1, \quad
\sqrt{2} \sum_{n=0}^\infty (-1)^n g_n \,=\, 0.\]
\end{lemma}

\begin{proof} 
In view of Corollary \ref{lemOrthnorm}, we have 
\begin{align}
1 \,=&\, \sum_{n=-\infty}^\infty |\widehat{\phi}(\lambda + 2n)|^2 \,=\,
\sum_{n=-\infty}^\infty |G(\lambda/2 + n)|^2\,
    |\widehat{\phi}(\lambda/2 + n)|^2\notag\\
 =&\,\, |G(\lambda/2)|^2\! \sum_{n=-\infty}^\infty 
  |\widehat{\phi}(\lambda/2 + 2n)|^2
\,+\, |G(\lambda/2 + 1)|^2 
\sum_{n=-\infty}^\infty |\widehat{\phi}(\lambda/2 + 2n + 1)|^2\notag\\
=& \,\, |G(\lambda/2)|^2 \,+\, |G(\lambda/2 + 1)|^2\notag
\end{align}
almost everywhere. If $\phi \in L^1(H)$ then $\,\widehat{\phi}\,$ is 
continuous and $\,\widehat{\phi}(0)\neq 0$ by Theorem \ref{thmMRA}.
 Hence $G(0)=1$ by (\ref{relG}) and  $G(1)=0$ is an immediate consequence of
 (\ref{Gsquare}).
\end{proof}

\section{Orthogonal radial wavelets}

In this section we construct wavelets for a given radial MRA $\{V_j\}_{j\in \Z}$ 
in $\R^3$ with
\textbf{orthonormal}  scaling function $\phi$ and filter function $G$.
As usual, the wavelet space $W_j$ is defined as the orthogonal complement of
$V_j$ in $V_{j+1}$, 
\[W_j:= V_{j+1} \ominus V_j\,.\]
Thus $L^2(H)$ decomposes as an orthogonal Hilbert sum 
\[L^2(H) \,=\, \bigoplus_{j=-\infty}^\infty W_j.\] 
Recall the definition of $S$ in Section 4  and the 
characterization 
of $V_0=V_\phi$ according to Lemma \ref{Lem_Contain},
\begin{equation}\label{contain2}
 f\in V_0 \,\Longleftrightarrow\,
\widehat f(\lambda)\,=\, \frac{\beta(\lambda)}{\lambda}\,
\widehat\phi(\lambda) \quad\text{with }\,\beta\in S.
\end{equation}
Define
\[ S_0:= \{\alpha\in S: \alpha(\lambda+1) = -\alpha(\lambda) \,\,\text{for almost all
}\, \lambda\}\]
which is a closed subspace of $S$ w.r.t. $\|\,.\|_{L^2[0,1]}$.
Then 
$W_{-1} = V_0\ominus V_{-1}$ is characterized as follows:

\begin{prop}\label{lem_W0cont} 
\begin{enumerate}\itemsep=-1pt
\item[\rm{(i)}] 
Let $f\in L^2(H).$ Then
\[ f\in W_{-1} \,\Longleftrightarrow\,\, 
\widehat{f}(\lambda) \,=\, \frac{\alpha(\lambda)}{\lambda}\, \ol{G(\lambda + 1)}\, 
\widehat{\phi}(\lambda) \,\text{ for some } 
\alpha\in S_0.\]
\item[\rm{(ii)}] The mapping $\, S_0 \to W_{-1}\,, \>\> \alpha \mapsto f_\alpha\,$ 
with 
 \[\widehat f_\alpha(\lambda):=
 (2\pi)^{1/4}\,\frac{\alpha(\lambda)}{\lambda}\, \ol{G(\lambda + 1)}\, 
\widehat{\phi}(\lambda)\]
is an isometric isomorphism.
\end{enumerate}
\end{prop}

\begin{proof} For (i), notice first that 
the Hankel transform is a unitary isomorphism of
$L^2(H),$ so $\,\widehat W_{-1} \,=\, \widehat V_{0}\ominus \widehat V_{-1}.$
 Rescaling of \eqref{contain2} by the factor $2$ and relation \eqref{relG} 
imply that $h\in L^2(H)$ is contained in 
$\widehat V_{-1}$  if and only if there  exists some $\widetilde\beta\in S$ 
such that
\[ h(\lambda) = \,\frac{\widetilde\beta(2\lambda)}{\lambda} \,\widehat\phi(2\lambda)
\,=\, 
\frac{\widetilde\beta(2\lambda)}{\lambda}\,G(\lambda)\,\widehat\phi(\lambda).\]
Thus $\beta\in S$ corresponds to $f\in W_{-1}$ according to \eqref{contain2}
if and only if
\[\int_0^\infty
\frac{\widetilde\beta(2\lambda)}{\lambda}\,G(\lambda)\,\frac{\ol{\beta(\lambda)}}{\lambda}
|\widehat\phi(\lambda)|^2 d\omega(\lambda) \,= \,0
   \quad \text{for all } \widetilde\beta\in S.\]
Up to a constant factor, the integral on the left equals
\begin{align}
\int_{-\infty}^\infty \widetilde\beta(2\lambda) G(\lambda)\,\ol{\beta(\lambda)}\, &
  |\widehat\phi(\lambda)|^2\,d\lambda \,=\,
 \sum_{n\in \Z}\int_0^1 \widetilde\beta(2\lambda)\,
\ol{\beta(\lambda+2n)}G(\lambda+2n)\,|\widehat\phi(\lambda+2n)|^2 d\lambda \notag\\
+\,&\, 
\sum_{n\in \Z}\int_0^1 \widetilde\beta(2\lambda)\, \ol{\beta(\lambda+2n+1)} 
G(\lambda+2n+1)\,|\widehat\phi(\lambda+2n+1)|^2 d\lambda \notag\\
=&\, \int_0^1 \widetilde\beta(2\lambda)\bigl(\ol{\beta(\lambda)}\,G(\lambda) +
\ol{\beta(\lambda +1)}
\,G(\lambda +1)\bigr) d\lambda,\notag
\end{align}
where we used the periodicity and symmetry properties of $\beta,\widetilde\beta, G$
as well
as Corollary \ref{lemOrthnorm}.
Since $\widetilde\beta\in S$ is arbitrary, we conclude that the vectors  
     $(\beta(\lambda),\beta(\lambda+1))^T$ and  $(G(\lambda),G(\lambda+1))^T$ must be
orthogonal in $\C^2$ for
almost all $\lambda$. This means that
\begin{equation}\label{orthorelation}
\left(\begin{array}{c} \beta(\lambda) \\ \beta(\lambda + 1) \end{array}\right)
\,=\, \alpha(\lambda) \left(\begin{array}{c} \ol{G(\lambda+ 1)} \\ -\ol{G(\lambda)}
\end{array}\right)
\end{equation}
for some  function $\alpha: [0,1]\to \C.$ Thanks to the boundedness of $G$,
$\alpha$  belongs to $L^2[0,1].$  Since $\beta$ and $G$ are $2$-periodic
on $\R$, $\beta$
is odd and $G$ is even, an extension of $\alpha$  to $\R$ 
 must be  $2$-periodic and odd; hence $\alpha\in S$ with 
$\beta(\lambda)= \alpha(\lambda)\ol{G(\lambda+1)}.$
Using the  $2$-periodicity of $\beta$ and $G$ in the second component of 
\eqref{orthorelation}, we further
deduce that $\beta(\lambda) = -\alpha(\lambda+1) \ol{G(\lambda +1)}$ and therefore
$\alpha(\lambda+1) = -\alpha(\lambda)$ apart from the zero-set of $G(\lambda +1)$.
Thus $\widehat f$ is of the claimed form. Conversely, if $\alpha\in S_0$,  then 
$\,\beta(\lambda):= \alpha(\lambda)\ol{G(\lambda+1)}\in S$, and 
\eqref{orthorelation} is satisfied.

\noindent
For the proof of (ii), we calculate
\begin{align} \|f_\alpha\|_2^2\,=\,\|\widehat f_\alpha\|_2^2\,=&\,
\sqrt{2\pi} \int_0^\infty 
  \frac{|\alpha(\lambda)|^2}{\lambda^2} \,|G(\lambda +1)|^2 
  |\widehat\phi(\lambda)|^2 d\omega(\lambda) \notag\\
 \,=&\, \int_{-\infty}^\infty |\alpha(\lambda)|^2|G(\lambda +1)|^2 
  |\widehat\phi(\lambda)|^2 d\lambda \,=\, \int_{-1}^1|\alpha(\lambda)|^2|G(\lambda
+1)|^2 d\lambda,\notag
\end{align}
where we used that $\alpha$ and $G$ are $2$-periodic and $\phi$ is an orthonormal
scaling 
function. 
By assumption on $\alpha$, we have $\alpha(\lambda-1) = \alpha(\lambda+1) =
-\alpha(\lambda)$. Thus by 
Lemma \ref{lemGprop},
\begin{align}
\int_{-1}^1|\alpha(\lambda)|^2|G(\lambda +1)|^2 d\lambda \,=&\, 
  \int_0^1 |\alpha(\lambda-1)|^2|G(\lambda)|^2 d\lambda \,+\, 
  \int_0^1 |\alpha(\lambda)|^2|G(\lambda +1)|^2 d\lambda\notag\\
  \,=&\, \int_0^1 |\alpha(\lambda)|^2d\lambda \,=\, \|\alpha\|_2^2\,.\notag
\end{align}
This proves (ii).
\end{proof}

It is now easy to obtain an orthonormal basis of $W_{-1}$. Recall
that the $s_k,\, k\in \N$ form an orthonormal basis of $S$. 
Moreover, the $s_{2k-1}, \, k\in \N$  are an orthonormal basis of $S_0$. 
Thus by the previous result, the functions 
\[f_k:= f_{s_{2k-1}} \quad(k\in \N)\]
constitute an orthonormal basis of $W_{-1}.$
Define $\psi\in L^2(H)$ by
\begin{equation}\label{def_wavelet}
\widehat{\psi}(2 \lambda)\,=\, \ol{G(\lambda +1)} \,\widehat \phi(\lambda).
\end{equation}
Then in view of \eqref{trans2},
\[ f_k \,=\,  \frac{M_{2k-1}}{2} T^{(2k-1)} D_2 \psi.\]
To obtain an orthonormal basis of $W_0$, we just have to rescale.  Extending the
notation $T^{(r)}:= T_{\pi r}$  to 
$\,r \in \N/2$ and using the relation $D_a T_x = T_{ax} D_a$, 
we obtain that an orthonormal basis of $W_0$ is given by the functions 
\[\psi_k:= D_{1/2} f_k =
\frac{M_{2k-1}}{2} T^{((2k-1)/2)}\psi, \quad k\in \N.\]
We call $\psi$ a (basic) wavelet for the radial mutiresolution
$(V_j).$

\begin{defn} For $j\in \Z$ and $k\in \N,$ define the ``radial'' wavelets
\[ \psi_{j,k}(r):= D_{2^{-j}} \psi_k(r) = 8^{j/2}  \frac{M_{2k-1}}{2}
T^{((2k-1)/2)}\psi (2^j r).\]
\end{defn}

 \noindent
We have proven:

\begin{thm}\label{thm_wavebasis} \begin{enumerate}\itemsep=-1pt
\item[\rm{(i)}] For each $j\in \Z$, the set
  $\,\{\psi_{j,k}:\, k\in \N\}\,$
   constitutes an orthonormal basis of $W_j$.
\item[\rm{(ii)}] The set $\{\psi_{j,k}:\, j\in \Z, \,k\in \N\}$ is
 an orthonormal wavelet basis of $L^2(H)$.
\end{enumerate}
\end{thm}

\begin{corollary} 
The functions
\[\Psi_{j,k} (x):= \psi_{j,k}(|x|),\quad x\in \R^3,\,j\in\Z,\, k\in \N\]
form an orthonormal basis for the closed subspace
$\, L_{rad}^2(\R^3) := \{f \in L^2(\R^3): \,f \text{ radial\,}\}\,$
of radial functions in  $L^2(\R^3)$. 
\end{corollary}

\section{Construction of radial scaling functions and wavelets}

Yet, we do not have handsome criteria in order to decide whether a given 
function $\phi\in L^2(H)$ is a radial scaling function, i.e. a scaling function for a
radial MRA. The analogy of
our constructions 
to those  on the group $(\R, +)$, however, leads to the following
close relationship:

\begin{thm}\label{thm_class} Suppose $\,\phi_\R$ is a classical scaling function 
on $\R$ which is even
 and such that its
(classical) Fourier transfrom $\F(\phi_\R)$
   is continuous in $0$ and satisfies 
 $\F(\phi_\R) \in L^2(H)$. 
Define  $\phi \in L^2(H)$ via its Hankel transform,
\begin{equation}\label{def_phi}
\widehat{\phi}(\lambda) := \sqrt{2\pi}\, \F(\phi_\R)(\pi \lambda).
\end{equation}
Then $\phi$ is a radial scaling function.

Conversely, if
$\phi$ is a scaling function for a radial MRA such that $\widehat{\phi}$
is continuous in $0$, then $\widehat \phi\in L^2(\R)$ and
 the function $\phi_\R$ defined by (\ref{def_phi}) 
(where $\widehat{\phi}$ is extended to
an even function on $\R$)
is a classical scaling function on $\R$. 

Moreover, $\phi$ is an \emph{orthonormal} radial scaling function 
if and only if $\phi_\R$ is an \emph{orthonormal} classical scaling function. 
\end{thm}

\begin{proof} Let us start with the first assertion. As $\phi_\R$ 
is a classical scaling function, we have 
by eq. (5.3.2) in \cite{Daub}
\[\frac{A}{2\pi} \leq \sum_{k\in \Z} |\F(\phi_\R)(\xi + 2\pi k)|^2 \leq
\frac{B}{2\pi} \quad \mbox{a.e.} \] 
with suitable constants $0<A\leq B\leq \infty$. Moreover, $\phi_\R$ is orthonormal
if and only if $A=B=1$. Since $\phi_\R$ is assumed to be even,  definition
(\ref{def_phi}) is compatible with the even extension of $\widehat{\phi}$.
By Proposition \ref{Prop_Riesz}, the set $\{M_k T^{(k)}\phi, \, k\in \N\}$ 
forms a Riesz
basis for $V_0 = \ol{\text{span} \{B_\phi\}}$ which is an orthonormal basis if and
only
if $\phi_\R$ is orthonormal; c.f. Corollary \ref{lemOrthnorm}. Moreover, by
eq. (5.3.18) of \cite{Daub}, there exists a $2\pi$-periodic function 
$m_0 \in L^2([-\pi,\pi])$ such that 
$\F(\phi_\R)(\xi) = m_0(\xi /2) \F(\phi_\R)(\xi / 2),$ and $m_0$ is  necessarily even
in our case. Hence, with 
$\gamma(\lambda) := m_0(\lambda \pi) \sin(2\pi \lambda)$ which clearly is
contained in $S,$ we have 
$\sin(2\pi \lambda) \widehat{\phi}(2\lambda) = \gamma(\lambda)
\widehat{\phi}(\lambda)$.
This is exactly the radial two-scale equation (\ref{two_scale}).
As $\widehat{\phi}$ is continuous in $0$, the condition $\widehat{\phi}(0) \neq 0$ 
of Theorem \ref{thmMRA} is automatically satisfied 
(see e.g. Remark 3 on p.144 in \cite{Daub}), and thus we finally obtain that $\phi$
is a radial scaling
function.

For the converse part notice first that continuity of $\,\widehat \phi\,$ in $0$
already implies  that $\widehat
\phi \in L^2(\R).$ We further proceed similar 
as before, using Proposition 5.3.1 and
Proposition 5.3.2 in \cite{Daub} and the corresponding results of the 
present paper. 
\end{proof}

This theorem supplies a variety of radial scaling functions since there are many
 classical scaling functions on $\R$ which  satisfy the assumptions
of the theorem. However, as to
orthonormal radial scaling functions with compact support, a 
famous theorem of Daubechies implies the following negative result.

\begin{corollary} There  do not exist any real-valued orthonormal 
radial scaling functions with compact support.
\end{corollary}

\begin{proof} The proof of Theorem 8.1.4 in \cite{Daub} shows that  an even, 
real-valued and
compactly supported scaling function  is
necessarily  the Haar function $\chi_{[-1/2,1/2]}$, 
the characteristic function of the 
interval $[-1/2,1/2]$. 
However, its Fourier transform
$\F(\chi_{[-1/2,1/2]})(\xi) = \sqrt{\frac{2}{\pi}}\frac{\sin(\xi/2)}{\xi}$
 is not contained in $L^2(H).$
\end{proof}

If $\phi$ correponds to an even classical scaling function $\phi_\R$ according
to Theorem \ref{thm_class}, then
 the hypergroup translates $\phi_{0,k}= M_k T_{k\pi}\phi\,$ may be expressed
accordingt to the formula
\begin{equation}\label{trans_phi}
\phi_{0,k}(x) \,=\, \frac{1}{(2\pi)^{1/4}\, x} \left(\phi_\R
\bigl(\frac{x}{\pi}-k\bigr) -
\phi_\R\bigl(\frac{x}{\pi} + k\bigr)\right). 
\end{equation}
In fact,  by the Plancherel theorem for the Hankel transform and 
eq. \eqref{TFB},
\[T_r \phi(s) = \,\int_0^\infty j(\lambda r) j(\lambda s) \widehat\phi(\lambda) 
  d\omega(\lambda)\,=\, \frac{1}{\sqrt{2\pi}\cdot rs} \int_{-\infty}^\infty
\widehat\phi(\lambda) \sin(s\lambda)
\sin(r\lambda)d\lambda \]
for all $r,s\in \R_+$. Here $\widehat\phi\,$ is as usual extended to
 an even function on $\R$. Using relation \eqref{def_phi} as well as basic
trigonometric identities and the
Plancherel theorem for the classical Fourier transform, 
 we can write
\begin{align}
 T_r \phi(s) \,=&\, \frac{1}{2\sqrt{2\pi}\cdot rs} \int_{-\infty}^\infty
\widehat\phi(\lambda) \bigl(\cos
\lambda(r-s) - \cos \lambda(r+s) \bigr)d\lambda \notag\\
=&\,\frac{1}{2rs}\int_{-\infty}^\infty \F(\phi_\R)(\pi\lambda)\left(e^{i\lambda(r-s)}
- e^{i\lambda(r+s)}
\right) d\lambda\notag\\ 
=&\,\frac{1}{\sqrt{2\pi}\cdot rs} \left(\phi_\R \Bigl(\frac{r-s}{\pi}\Bigr) - 
\phi_\R \Bigl(\frac{r+s}{\pi}\Bigr)\right)\notag
\end{align}
This implies \eqref{trans_phi}.

As an example, we consider the radial analogue of the Shannon wavelets.
We define the scaling function via its Hankel transform,
\[\widehat{\phi}(\lambda) \,=\, \chi_{[0,1]}(\lambda),\qquad
 \phi(x) %= \sqrt{\frac{2}{9\pi}} j_{3/2}(x) 
 =\,  \sqrt{\frac{2}{\pi}}  \frac{\sin(x) - x\cos(x)}{x^3}.\]
Constructing the associated basic wavelet according to formula (\ref{def_wavelet})
yields (after
a short calculation)
\[\hat{\psi}(\lambda) =\, \chi_{[1,2]}(\lambda),\qquad
\psi(x) \,=\, \sqrt{\frac{2}{\pi}}\frac{\sin(2x)-\sin(x) -2x\cos(2x) +
x\cos(x)}{x^3}.
\]
The translates of the scaling function and the wavelet turn out to be
\begin{align}
\phi_{0,k}(x) =&\, \frac{1}{(2\pi)^{1/4} x}\left(\frac{\sin(x-k\pi)}{x-k\pi} - 
\frac{\sin(x+k\pi)}{x+k\pi}\right),~k\in\N,\notag\\
\psi_{0,k}(x) =&\, \frac{1}{(2\pi)^{1/4} x}\left(p\left(x - \frac{2k-1}{2}\pi\right)
- p\left(x+\frac{2k-1}{2}\pi\right)\right),~k\in\N,\notag
\end{align}
with $\,\displaystyle p(x)= \frac{\sin(2x)-\sin(x)}{x}$.

\section{Algorithms}

For the use of our radial multiresolution in applications
we need to formulate decomposition and
reconstruction algorithms. 
The first step in such an algorithm consists of projecting the 
function $f$ into a scale space $V_j$ for some suitable $j$. We obtain a
representation
\begin{equation*}%\label{projf}
 P_j f \,=\, \sum_{k=1}^\infty c_k^{(j)} \phi_{j,k} .
\end{equation*}
So from now on we assume that we have given a function $f\in V_j$ 
in terms of its coefficients $c_k^{(j)}$. 
The decomposition algorithm consists of
decomposing $f$ into $V_{j-1}$ and $W_{j-1}$, i.e. of calculating the
coefficients $c_k^{(j-1)}$ and $d_k^{(j-1)}$ in the representation
\begin{equation*}\label{repf}
f ~=~ \sum_{k=1}^\infty c_k^{(j-1)} \phi_{j-1,k} + 
\sum_{k=1}^\infty d_k^{(j-1)} \psi_{j-1,k}.
\end{equation*}
(Such a representation exists, since by construction $\{\phi_{j-1,k},\psi_{j-1,k}~:~k\in\N\}$ is also a basis of $V_j$.) 
A reconstruction
algorithm determines the coefficients $c_k^{(j)}$ when
$f$ is given in terms of $c_k^{(j-1)}$ and $d_k^{(j-1)}$, $k\in \N$.

We still assume that $\phi$ is an orthonormal scaling function (and $\psi$
is hence an orthonormal wavelet). Let
\begin{align}
q_\ell^{(k)} \,:=\, \langle \phi_{1,k}, \phi_{0,l} \rangle \,=\, 
\langle \phi_{j,k}, \phi_{j-1,l} \rangle,\qquad
r_\ell^{(k)} \,:=\, \langle \phi_{1,k}, \psi_{0,l} \rangle \,=\,
\langle \phi_{j,k}, \psi_{j-1,l} \rangle.\nonumber
\end{align}
By using Hilbert space techniques - in particular
Parseval's equation - we obtain analogously as in standard wavelet theory 
the decomposition formulae
\begin{align}
c_\ell^{(j-1)} \,=\, \sum_{k=1}^\infty c_k^{(j)} q_\ell^{(k)},\qquad
d_\ell^{(j-1)} \,=\, \sum_{k=1}^\infty c_k^{(j)} r_\ell^{(k)},\notag
\end{align}
and the reconstruction formula
\begin{equation*}
c_k^{(j)} ~=~ \sum_{\ell=1}^\infty c_\ell^{(j-1)} 
\overline{q_\ell^{(k)}} + \sum_{\ell=1}^\infty d_\ell^{(j-1)}
\overline{r_\ell^{(k)}}.  
\end{equation*}
It turns out that the coefficients $q_\ell^{(k)}$ and
$r_\ell^{(k)}$ are determined in terms of the numbers $g_n$ in the cosine expansion
of $G$, i.e. the coefficients in
\[G(\lambda) ~=~ \sqrt{2} \sum_{n=0}^\infty g_n \cos(n\pi \lambda).\]

\begin{thm} For $\ell,k \in \N$ it holds
\begin{align}
q_\ell^{(k)} \,=\,& \left\{\begin{array}{ll}\displaystyle 
\overline{g_{k-2\ell} - g_{2\ell + k}} & \mbox{\rm for } 
2\ell < k,\\
\overline{2g_0 - g_{4\ell}} & \mbox{\rm for } 2\ell = k, \\
\overline{g_{2\ell-k} - g_{2\ell + k}} & \mbox{\rm for } 2\ell > k,
\end{array}\right. \nonumber\\
r_\ell^{(k)} \,=\,&  \left\{\begin{array}{ll}
 (-1)^{k-1}(g_{k-2\ell+1} - g_{k+2\ell-1}) & \mbox{\rm for } 2\ell - 1 < k,\\
 2 g_0 - g_{4\ell-2} & \mbox{\rm for } 2\ell-1 = k, \\
 (-1)^{k-1}(g_{2\ell-1-k} - g_{2\ell - 1 + k}) & \mbox{\rm for } 2\ell - 1 > k.
\end{array}\right.\nonumber
\end{align}
\end{thm}

\begin{proof} Using the Plancherel theorem, relation (\ref{relG}) and Corollary
(\ref{lemOrthnorm}) 
we obtain
\begin{align}
q_\ell^{(j)} \,=\,& \langle \phi_{1,k}, \phi_{0,\ell} \rangle
\,=\, \langle \hat{\phi}_{1,k}, \hat{\phi}_{0,\ell} \rangle\nonumber\\
\,=\,& 8^{-1/2} \int_0^\infty \rho_k(\lambda/2) \rho_\ell(\lambda/2)
\overline{G(\lambda/2)}
|\hat{\phi}(\lambda/2)|^2 d\omega(\lambda)\nonumber\\
\,=\,& \sqrt{2} \int_0^1 s_k(\lambda) s_\ell(2\lambda) \overline{G(\lambda)} 
\sum_{n=-\infty}^\infty |\hat{\phi}(\lambda + 2n)|^2 d\lambda\nonumber\\
\label{repq}
\,=\,& 4 \sum_{n=0}^\infty \ol{g_n} \int_0^1 
\sin (k\lambda \pi) \sin(2\ell \lambda \pi) \cos(n\pi \lambda) d\lambda.
\end{align}
An easy calculation using trigonometric identities shows
\begin{equation}\label{intcalc}
\int_0^1 \sin (k\lambda \pi) \sin(t \lambda \pi) \cos(n\pi \lambda) d\lambda
\,=\, \left\{\begin{array}{ll} \frac{1}{4}(\delta_{n,|t-k|} - 
\delta_{n,t+k}) & \mbox{for } n> 0,\\
\frac{1}{4}(2\delta_{0,t-k} - \delta_{0,t + k}) & \mbox{for } n=0.
\end{array}\right.
\end{equation}
Setting $t=2\ell$ and inserting into (\ref{repq}) yields 
the assertion for $q_\ell^{(k)}$. We proceed similarly for
$r_\ell^{(k)}$:
\begin{align}
r_\ell^{(j)} \,=\,& \langle \phi_{1,k}, \psi_{0,\ell} \rangle
\,=\, \sqrt{2} \int_0^1 s_k(\lambda) s_{2\ell-1}(\lambda) G(\lambda+1) 
d\lambda\nonumber\\
\,=\,& 4 \sum_{n=0}^\infty g_n \int_0^1 \sin(k\pi\lambda) 
\sin((2\ell-1)\pi\lambda) \cos(n\pi(\lambda + 1)) d\lambda\nonumber\\
\,=\,& 4 \sum_{n=0}^\infty g_n (-1)^n \int_0^1 \sin(k\pi\lambda) 
\sin((2\ell-1)\pi\lambda) \cos(n\pi\lambda) d\lambda.\nonumber
\end{align}
Setting $t=2\ell -1$ in (\ref{intcalc}) and inserting into the 
last expression gives the result for $r_\ell^{(k)}$.
\end{proof}

Let us consider the case where only finitely many 
coefficients $g_k$ are different from zero. Although this 
is not possible for real-valued
orthonormal scaling functions this assumption makes it easier to compare the
radial wavelet algorithm with the classical one.
Of course, in applications one can only handle
finitely many coefficients anyway.
So let us assume $\supp g \subset [0, N]$, i.e. $g_k =0$ for 
$k \notin \{0,\hdots,N\}$. Elementary considerations show the following.
Leaving $k$ fixed yields
\begin{equation*}
\begin{array}{lclll}
q_\ell^{(k)} &=& 0 &\mbox{for } \ell 
\notin [\frac{k-N}{2},\frac{k+N}{2}] &
\mbox{if } k > N,\\
q_\ell^{(k)} &=& 0 & \mbox{for }  \ell \notin [1,\frac{k+N}{2}] &
\mbox{if } k \leq N,\\
r_\ell^{(k)} &=& 0 & \mbox{for } \ell \notin [\frac{k-N-1}{2},
\frac{N+k+1}{2}] & \mbox{if }k > N+1,\\
r_\ell^{(k)} &=& 0 & \mbox{for } \ell \notin [1,\frac{N+k+1}{2}] &
\mbox{if }k \leq N+1.
\end{array}
\end{equation*}
If $\ell$ is fixed then
\begin{equation*}
\begin{array}{lclll}
q_\ell^{(k)} &=& 0 & \mbox{for } k \notin [2\ell - N,2\ell+N] &
\mbox{if } 2\ell > N,\\
q_\ell^{(k)} &=& 0 & \mbox{for } k \notin [1,N+2\ell] &
\mbox{if } 2\ell \leq N,\\
r_\ell^{(k)} &=& 0 & \mbox{for } k \notin [2\ell-1-N,2\ell-1+N] &
\mbox{if } 2\ell -1 > N,\\
r_\ell^{(k)} &=& 0 & \mbox{for } k \notin [1,2\ell-1+N] &
\mbox{if } 2\ell -1 \leq N.
\end{array}
\end{equation*}
With
\begin{equation*}
h_k \,:=\, \left\{\begin{array}{ll} 
g_{|k|} & \mbox{for } 1\leq |k| \leq N-1,\\
2g_0 & \mbox{for } k=0,\\
0 & \mbox{otherwise}.
\end{array}\right.
\end{equation*}
it holds
$G(\lambda) ~=~ \frac{1}{\sqrt{2}} \sum_{k=-N}^N h_k e^{ik\lambda}$.
Because of the conditions on it, $G$ is also the filter function
for an ordinary multiresolution analysis on $\R$ with coefficients $h_k$.
Now, if $2\ell > N$ resp. $2\ell -1 > N$ then it is easy to see that
\begin{align}
q_\ell^{(2\ell+k)}   \,=\,& \ol{h_k} \quad \mbox{for } k=-N,\hdots,N,\nonumber\\
r_\ell^{(2\ell-1+k)} \,=\,& (-1)^k h_k \quad \mbox{for } k=-N,\hdots,N.\nonumber
\end{align}
Similarly, if $k > N$ then
\begin{align}
q_\ell^{(k)} \,=\, \ol{h_{k-2\ell}},\qquad
r_\ell^{(k)} \,=\,  (-1)^{k+1} h_{k+1-2\ell}.\nonumber 
\end{align}
Hence, for $2\ell > N + 1$ the decomposition formulae become
\begin{equation*}
c_\ell^{(j-1)} \,=\, \sum_{k=1}^\infty c_k^{(j)} \ol{h_{k-2\ell}},\qquad 
d_\ell^{(j-1)} \,=\, \sum_{k=1}^\infty c_k^{(j)} (-1)^{k+1} h_{k+1-2\ell},
\end{equation*}
and for $k > N+1$ the reconstruction formula is 
\begin{equation*}
c_k^{(j)} \,=\, \sum_{\ell=1}^\infty c_\ell^{(j-1)} h_{k-2\ell} 
+ \sum_{\ell=1}^\infty d_\ell^{(j-1)} (-1)^{k+1} \ol{h_{k+1-2\ell}}.
\end{equation*}
These formulae are well-known. Indeed, they are the decomposition and reconstruction
formulae of the classical discrete wavelet transform. 
So our approach leads to the classical algorithm if we are
far enough away from the origin. If we are close to the origin we have
derived an algorithm to handle the boundary point $0$.

\noindent
Holger Rauhut\\
Zentrum Mathematik\\
Technische Universit\"at M\"unchen\\
D--80290 M\"unchen\\
Germany\\
rauhut@ma.tum.de\\

\medskip\noindent
Margit R\"osler\\
Mathematisches Institut\\
Universit\"at G\"ottingen\\
Bunsenstr. 3-5\\
D-37073 G\"ottingen\\
Germany\\ 
roesler@uni-math.gwdg.de

\end{document}